\documentclass[12pt]{amsart}
\usepackage{float}
\usepackage{booktabs}
\usepackage{longtable}

\setlongtables

\usepackage{fancyhdr}
\newtheorem{defn}{Definition}
\newtheorem{conj}{Conjecture}
\newtheorem{question}{Question}

\newcommand{\up}{\uparrow}

\newcommand{\rarrowf}[1]{{\hbox to #1{\rightarrowfill}}}

\newcommand{\ov}{\overline}

\newcommand{\Q}{\ensuremath{\mathbb Q}}

\DeclareTextSymbol{\ij}{T1}{188}
\DeclareTextCommandDefault{\ij}{ij}

\newcommand{\C}{{\mathbb C}}

\newcommand{\tensor}{\otimes}

\newcommand{\Qp}{\ensuremath {{\mathbb Q}_{p}}}
\newcommand{\Cp}{\ensuremath {{\mathbb C}_{p}}}

\newcommand{\padic}{$p$-adic}

\newcommand{\Tp}{\ensuremath{\mathrm{T}_{\!p}}}

\DeclareMathOperator{\ssing}{ss}
\DeclareMathOperator{\slope}{slope}

\DeclareMathOperator{\ord}{ord}

\newcommand{\U}{\mathrm{U}}

\def\m@th{\mathsurround=0pt }
\def\n@space{\nulldelimiterspace=0pt \m@th}

\def\mDown#1#2{{\hbox{$\left\downarrow\vbox to #1{}\right.\n@space$}}
         \rlap{$\vcenter{\hbox{$\scriptstyle#2$}}$}}
\def\mUp#1#2{{\hbox{$\left\up\vbox to #1{}\right.\n@space$}}
         \rlap{$\vcenter{\hbox{$\scriptstyle#2$}}$}}
\def\meq#1{{\hbox{$\left\|\vbox to #1{}\right.\n@space$}}}

\DeclareMathAlphabet{\mathbb}{U}{msb}{m}{n}

\newcommand{\mytoday}{\ifcase\month\or
 january\or february\or march\or april\or may\or june\or
 july\or august\or september\or october\or november\or december\fi
 \space{\number\day}, {\number\year}}

\begin{document}

\title{Where the Slopes Are}
\author{Fernando Q. Gouv\^ea}
\address{Department of Mathematics and
Computer Science\\Colby College\\Waterville, ME 04901\\U. S. A.}
\email{fqgouvea@colby.edu}
\thanks{Project sponsored in part by the National Security Agency under
  Grant number MDA904-98-1-0012. The United States Government is authorized
  to reproduce and distribute reprints.}

\maketitle

Let $N$ be a positive integer (the ``level'') and $k\geq 2$ be an integer
(the ``weight''). Let $S_k(N,\C)$ denote the finite-dimensional $\C$-vector
space of cuspidal modular forms of weight $k$ and trivial character on
$\Gamma_0(N)$ defined over $\C$. Elements $f\in S_k(N,\C)$ can be specified
by giving their Fourier expansions
\[ f = a_1q + a_2q^2 + \dots = \sum_{n=0}^\infty a_nq^n,\]
where $q=e^{2\pi i z}$ and $z$ is in the complex upper halfplane. This
expansion is sometimes described as ``the $q$-expansion at infinity'' of
the modular form $f$. There exists a natural basis of $S_k(N,\C)$
consisting of forms all of whose Fourier coefficients are in fact
rational. We denote the $\Q$-vector space spanned by this basis by
$S_k(N,\Q)$. Note that then we have \[ S_k(N,\C) = S_k(N,\Q) \tensor \C.\] 

For each prime number $p$ which does not divide $N$ there is a linear
operator $\Tp$ acting on $S_k(N,\C)$, known as the $p$-th Hecke operator.
(In fact, the $\Tp$ stabilize $S_k(N,\Q)$.) A modular form which is an
eigenvector for all of these linear operators simultaneously is called an
\emph{eigenform}; the space $S_k(N,\C)$ has a basis made up of eigenforms,
and the Fourier coefficients of these eigenforms can be normalized (by
requiring $a_1=1$) to belong to a finite extension of $\Q$. 

The eigenvalues of the $\Tp$ operator encode significant arithmetic
information about the modular form and various other objects which can be
attached to it (for example, a Galois representation). In our setting, the
eigenvalue of $\Tp$ acting on an eigenform $f\in S_k(N,\Cp)$ is a totally
real algebraic number whose absolute value (with respect to any embedding
of $\ov\Q$ into $\C$) is between $-2p^{(k-1)/2}$ and $2p^{(k-1)/2}$. If we
normalize the eigenvalues by dividing by $p^{(k-1)/2}$, the normalized
eigenvalues are real numbers in the interval $[-2,2]$, and we can ask about
their distribution in that interval. The Sato-Tate Conjecture, still very
much an open problem, predicts the properties of that distribution for a
fixed modular form $f$ and varying $p$. We can also, however, fix the prime
$p$ and consider the distribution as $k\to\infty$ of all the eigenvalues of
$\Tp$ corresponding to eigenforms of weight $k$. This was done by Serre in
\cite{serre-distrib} and by Conrey, Duke, and Farmer in \cite{cdf}.

The goal of this paper is to begin the study of an analogous question in
the $p$-adic setting by presenting a wide range of numerical data. The
unexpected regularities in the data suggest several interesting questions
that deserve further investigation.

We fix a prime number $p$, then, and consider the situation in a $p$-adic
setting. We choose an embedding of the algebraic closure of $\Q$ into the
completion $\Cp$ of an algebraic closure of $\Qp$, and then we define
\[ S_k(N,\Cp) = S_k(N,\Q) \tensor \Cp,\]
and similarly for $S_k(N,F)$ where $F$ is any extension of $\Qp$. In the
$p$-adic context, it turns out that the right operator to consider is not
$\Tp$ but rather the Atkin-Lehner $U$ operator, which can be described by 
its action on $q$-expansions:
\[ U\left(\sum a_nq^n\right) = \sum a_{np} q^n.\]
If $p$ does not divide $N$, this operator does not stabilize the space
$S_k(N,F)$, but it does stabilize the larger space $S_k(Np,F)$, and once
again we can consider eigenforms and the corresponding eigenvalues of $U$. 

Let $f\in S_k(Np,\Cp)$ be an eigenform for $U$, so that $U(f)=\lambda f$.
The $p$-adic valuation of the eigenvalue $\lambda$ turns out to play a
crucial role in the $p$-adic theory. We shall call this valuation the
\emph{slope} of the eigenform $f$:
\begin{defn} Given an $U$-eigenform $f$ of level $Np$, weight $k$ and
  eigenvalue $\lambda$, we define the \emph{slope} of $f$ by
  \[ \slope(f)=\ord_p(\lambda).\]
\end{defn}
The name ``slope'' comes from the $p$-adic theory of Newton polygons: the
slopes of the eigenforms in $S_k(Np,\Cp)$ are determined by the slopes of
the Newton polygon of the characteristic polynomial of the $U$ operator
acting on this space.

We are interested in the distribution \emph{of the slopes} of the $U$
operator for fixed level and varying weight. (Thus, we are writing the
eigenvalues as a $p$-adic unit times a power of $p$, and then we are
ignoring the unit part.) All of our results are numerical, but we feel they
are of sufficient interest and that they raise significant questions that
need to be addressed on a theoretical level.

I am grateful to several people for their contributions to this work.  The
main question discussed in this paper was raised by Dipendra Prasad in
conversation with the author. The computations were done with the
\textsc{gp} program using a modified version of a script written by Robert
Coleman. Finally, Barry Mazur, Kevin Buzzard, and Naomi Jochnowitz made
significant suggestions and observations at several points.

\section{Setting up the Problem}

Let $p$ be a prime number, $k\geq 2$ an even integer, and $N$ a positive
integer not divisible by $p$. Let $\ord_p$ be the $p$-adic valuation
mapping, normalized by $\ord_p(p)=1$. For any field $F$ of characteristic
zero, we write $S_k(N,F)$ to denote the $F$-vector space of cuspidal
modular forms of weight $k$ for $\Gamma_0(N)$ (with trivial character)
whose Fourier coefficients all belong to $F$. We will essentially be
concerned only with $F=\Qp$, since the Newton polygon (and therefore the
slopes) can be computed already in this context, though the eigenforms
themselves may only be defined over some extension of $\Qp$. Our
computations will be restricted to the case $N=1$ (and hence $k\geq 12$),
but it seems reasonable to set up the problem for general level.

There are two natural inclusions of $S_k(N,F)$ into $S_k(Np,F)$; on
$q$-expansions the first is the identity mapping and the second is the
Atkin-Lehner $V$ operator, which sends $q$ to $q^p$. The subspace spanned
by the images of both maps is called the space of \emph{oldforms} in
$S_k(Np,F)$; it has a natural complement called the space of
\emph{newforms}.

The Atkin-Lehner $U$ operator maps $S_k(Np,F)$ to itself, acting on
$q$-expansions by
\[ U(\sum a_nq^n) = \sum a_{np} q^n.\]
It follows from the Atkin-Lehner theory of change of level (see
\cite{atkin-lehner}) that the action of $U$ on newforms can be diagonalized
(possibly after extending the base field), and that all the eigenvalues are
equal to $\pm p^{(k-2)/2}$, and hence have slope equal to $(k-2)/2$. Thus,
as far as the slopes are concerned, the interesting questions have to do
with the action of $U$ on the oldforms. This is best understood by relating
it to the action of the Hecke operator $\Tp$ on forms of level $N$; this
yields the theory of ``twin eigenforms'' discussed in \cite{conj}.

The Hecke operator $\Tp$ can be diagonalized on $S_k(N,\Cp)$. Let $f\in
S_k(N,\Cp)$ be a normalized cuspidal eigenform, and let $a_p$ be the
eigenvalue of $\Tp$ acting on $f$. Finally, let $f_1, f_2\in S_k(Np,\Cp)$
be the two images of $f$ under the maps described above. The $U$ operator
stabilizes the two-dimensional space generated by $f_1$ and $f_2$, and its
characteristic polynomial is $x^2-a_p x + p^{k-1}$. If this polynomial has
two distinct roots, the action of $U$ on this two-dimensional subspace can
be diagonalized, and the slopes of the two resulting eigenforms can be
easily determined:
\begin{itemize}
\item If $\ord_p(a_p)<(k-1)/2$, the two eigenvalues have $p$-adic valuation
  equal to $\ord_p(a_p)$ and $k-1-\ord_p(a_p)$.
\item If $\ord_p(a_p)\geq (k-1)/2$, then both eigenvalues have $p$-adic
  valuation $(k-1)/2$.
\end{itemize}

It has been conjectured by Ulmer that the polynomial $x^2-a_px+p^{k-1}$
always has two distinct roots. Specifically:

\begin{conj}[Ulmer] The action of $U_p$ on $S_k(\Gamma_0(Np),\Q_p)$
  is semi\-sim\-ple. In particular, the polynomial $x^2-a_p x + p^{k-1}$
  always has distinct roots.
\end{conj}

Coleman and Edixhoven have shown that this is true for $k=2$ and that for
general $k$ it follows from the Tate Conjecture (see
\cite{coleman-edixhoven}). 

It is easy to see that if the polynomial has a double root then we must
have $\ord(a_p)=(k-1)/2$ (unless $p=2$, in which case we would have
$\ord(a_p)=(k+1)/2$).  In our computations (which were restricted to the
case $N=1$), \emph{we never found a form for which $\ord(a_p)\geq
  (k-1)/2$}, so that Ulmer's conjecture is true in every case we
considered. This shows that the computations below are correct independent
of Ulmer's conjecture.

For what follows, \emph{we assume Ulmer's conjecture holds}. In that case,
we can determine the slopes of $U$ acting on the oldforms in $S_k(Np,\Qp)$
by determining the slopes of $\Tp$ acting on $S_k(N,\Qp)$. For each slope
obtained in level $N$ one obtains a pair of slopes $\alpha'$ and $\alpha''$
in level $Np$, satisfying
\begin{itemize}
\item $0\leq \alpha'\leq \alpha'' \leq k-1$
\item $\alpha' + \alpha'' = k-1$
\end{itemize}
with $\alpha' < \alpha''$ unless they are both equal to $(k-1)/2$. 
We define the \emph{slope sequence for level $N$, weight $k$, and prime
  $p$} to be the ordered list of slopes
\[ (\alpha_1, \alpha_2, \dots, (k-1)-\alpha_2, (k-1)-\alpha_1) \]
for $U$ acting on the oldforms in $S_k(Np,\Qp)$, where we repeat slopes
that occur with multiplicity. The number of elements in this sequence is
equal to twice the dimension of $S_k(N,\Qp)$. Since the slope sequence is
symmetric under $\alpha \leftrightarrow (k-1)-\alpha$, we will usually
specify it by giving only the first half of the slope sequence.  The
discussion above shows that this first half is the same as the slope
sequence for $\Tp$ acting on $S_k(N,\Qp)$, except that all slopes larger
than $(k-1)/2$ are replaced by $(k-1)/2$. (As mentioned above, we found no
example in which the latter case occurs.)

Since we know that the slopes are in the interval $[0,k-1]$ (and we want to
vary $k$), it makes sense to normalize the slopes by dividing them by
$k-1$.

\begin{defn} Suppose $f$ is either a $\Tp$-eigenform of level $N$ or a
  $U$-eigenform $f$ of level $Np$. Let $k$ be the weight of $f$ and let
  $a_p(f)$ be the eigenvalue (of $\Tp$ or of $U$). We define the
  \emph{supersingularity} of $f$ by 
  \[ \ssing(f)=\frac{\ord_p(a_p(f))}{k-1}.\]
\end{defn}

Let $f\in S_k(N,\Cp)$ is an eigenform for $\Tp$ and (still
assuming Ulmer's Conjecture is true) $f',f''$ are the two old
$u$-eigenforms corresponding to it as above. Then, provided that 
$\ssing(f)\leq 1/2$, we have
\[ \ssing(f') = \ssing(f) \]
and
\[ \ssing(f'') = 1 - \ssing(f),\]
and both numbers are in the interval $[0,1]$. Thus, the sequence of
supersingularities corresponding to old eigenforms of weight $k$ and level
$Np$ is a normalized version of the slope sequence, and can be computed via
the supersingularities of forms of level $N$, provided these are small
enough. 

(One can think of $\ssing$ as a function on the \emph{eigencurve} studied
by Coleman and Mazur in \cite{eigencurve}. It will be an continuous
function on the eigencurve, except along the $k=1$ locus. Notice, however,
that classical eigenforms of weight $1$ will always have slope zero;
defining $\ssing(f)=0$ for such forms gives a continuous extension of
$\ssing$ to classical forms of weight $1$. No such continuous extension is
possible at points corresponding to non-ordinary forms of weight $1$.)

We define the \emph{supersingularity sequence in weight $k$} 
\[ (\eta_1, \eta_2, \dots, 1-\eta_2, 1-\eta_1) \]
by
\[ \eta_i = \ssing(f_i) = \frac{\slope(f_i)}{k-1}\]
as $f_i$ runs through the old eigenforms of weight $k$ on $\Gamma_0(Np)$.
The supersingularity sequence is contained in the interval $[0,1]$ and is
symmetric under $\eta \leftrightarrow 1-\eta$.

The main problem we consider is to understand the distribution of the
supersingularities in the interval $[0,1]$ when we fix the level $N$ and let
$k\to\infty$. This problem can be expressed in measure-theoretic terms, as
in \cite{serre-distrib}: considering $N$ and $p$ as fixed, for each $k$ we
define a probability measure $\mu_k$ on the interval $[0,1]$ by putting a
point mass at each supersingularity $\eta_i$: let $d_k=\dim S_k(N,\Qp)$,
and set
\[ \mu_k = \frac{1}{2d_k} \sum_{i=1}^{d_k} 
\left( \delta_{\eta_i} + \delta_{1-\eta_i}\right) , \] where $\delta_{x}$
is the Dirac measure at $x$. The question then is whether the measures
$\mu_k$ tend to a limit as $k\to\infty$, and if so to determine that limit
measure.

\section{Computations}

For our computations, we restricted to the case $N=1$, which then means
that one only gets non-trivial results for even weights $k\geq 12$. For
each prime number $p\leq 100$, we computed the Newton polygon of $\Tp$
acting on forms of weight $k$ and level $1$ for weights $k\leq 500$.  Since
in every case the slopes were less than $(k-1)/2$, the slopes we obtained
are exactly the first half of the slope sequence for the $U$ operator
acting on oldforms of level $p$, as described above.

The method used for computation was straightforward: the space of modular
forms of weight $k$ and level $1$ has a basis consisting of forms
$E_4^aE_6^b\Delta$, where $E_4$ and $E_6$ are the Eisenstein series of
weight $4$ and $6$ respectively, $\Delta$ is the unique cuspform of weight
$12$, and $4a+6b+12=k$. Using this explicit basis we determined the
characteristic polynomial of $\Tp$ and computed its Newton slopes, then
produced supersingularities by dividing by $k-1$. The computation was done
with the GP calculator \cite{gp-pari}; the basic GP functions we needed
were based on a script originally written by Robert Coleman. The main
constraint on the computation was the memory required for computing the
characteristic polynomial: larger $k$ meant working with a larger basis,
and larger $p$ meant that we needed to use more terms from the
$q$-expansion of the modular forms.  The full output of the computations
can be found on the web at
\texttt{http://www.colby.edu/personal/fqgouvea/slopes/}.

As already mentioned above, in every case we found that every slope in the
Newton polygon of $\Tp$ acting on forms of level $N$ was smaller than
$(k-1)/2$, from which it follows that $U$ acts semisimply on the space of
oldforms and that the slope sequence for $\Tp$ is indeed the same as the
first half of the slope sequence of $U$ acting on oldforms of level $Np$.

\begin{question}
  Fix a prime number $p$ and a level $N$. Let $f\in S_k(N,\Cp)$ be an
  eigenform for $\Tp$. Is it true that
  \[ \slope(f_i) < \frac{k-1}{2} \]
  always?
\end{question}

In fact, one sees much more. Even a cursory observation of the tables
suggests that the slopes are much smaller than one might expect. In fact,
we found that \emph{in almost every case} the supersingularities for weight
$k$ and prime $p$ are smaller than $1/(p+1)$. In other words, the
inequality
\[ \ssing(f) \leq \frac{1}{p+1} \]
holds almost always for forms of level $N$. It follows that the sequence of
supersingularities for weight $k$ is almost always contained in
$[0,\frac{1}{p+1}]\cup[\frac{p}{p+1},1]$. 

\begin{question}
  Fix a prime number $p$ and a level $N$.  Let $f\in S_k(N,\Cp)$ be an
  eigenform for $\Tp$. Is it true that
  \[ \ssing(f) \leq \frac{1}{p+1} \]
  almost always as $k\to\infty$?
\end{question}

\begin{table}

\begin{tabular}{cl}
\toprule
Prime $p$ & Weights $k$ \\ \midrule
59         & 16, 46, 76, 106, 136, 166, 196, 226, 256\\
           & 286, 316, 346, 376, 406, 436, 466, 496 \\
79         & 38, 44, 118, 124, 198, 204, 278, 284\\
           & 358, 364, 438, 444 \\
2411       & 12\\
15271      & 16\\
187441     & 16\\
3371       & 20\\
64709      & 20\\
27310421   & 26\\
\bottomrule
\end{tabular}

\caption{Known exceptions to $\eta_i\leq 1/(p+1)$}\label{exceptions}
\end{table}

To be more explicit about ``almost always,'' in our computations exceptions
to this inequality occurred only for $p=59$ and $p=79$; for each of these
primes, the inequality fails to hold for the highest-slope form in certain
weights. See Table~\ref{exceptions} for the list of weights at which
exceptional slopes appear; we discuss this list of weights further below.
Other exceptions to the inequality, outside the range of this computation,
can be read off from the results in \cite{story}; they correspond to forms
of weights $k=12, 16, 20$ that are non-ordinary with respect to large
primes. The final entry in the table comes from a computation by Atkin. The
full list of primes and weights for which we know of a slope that does not
satisfy the inequality is given in Table~\ref{exceptions}. For each $(p,k)$
pair, we found that exactly one slope in the first half of the slope
sequence violates the inequality.

If we focus on the exceptional slopes and compute the corresponding
supersingularities, we see that in our examples $\ssing(f)$ seems to get
closer to $1/(p+1)$ as $k$ grows. Let $p=59$, for example; the sequence of
supersingularities corresponding to the exceptional slopes in
Table~\ref{exceptions} is
\begin{align*}
&0.066, 0.022, 0.026, 0.019, 0.022, 0.018, 0.020, 0.017, 0.019, \\
&0.017, 0.019, 0.017, 0.018, 0.017, 0.018, 0.017, 0.018
\end{align*}
Here, of course, $1/(p+1)=1/60=0.01666\dots$, and the exceptional values of
the supersingularity seem to be approaching this value as the weight grows.

Encouraged by this, we can try to make ``almost always'' precise by using
the measure-theoretic formulation:

\begin{question}
  Is it true that the sequence of measures $\{\mu_k\}$ converges, as
  $k\to\infty$, to a measure supported on the set $[0,\frac{1}{p+1}] \cup
  [\frac{p}{p+1},1]$?
\end{question}

One way to think about these results is by analogy with the results in
\cite{story}. There, we fixed a modular form and computed its slopes with
respect to varying primes $p$, and found that the slope was almost always
zero. In this case, we fix a prime and consider the full slope sequences
for varying weights $k$, and we find that the slopes are almost always
bounded by $(k-1)/(p+1)$. The two exceptional sets are connected, of
course. For example, a form of weight $k$ which has slope equal to $1$ at a
prime $p>k+2$ will be exceptional from both points of view.

We know of no general results which suggest that ``exceptional'' slopes are
rare. For specific primes, there are some hints. For $p=3$ and $k=2\cdot
3^a$, Lawren Smithline has shown in \cite{smithline-thesis} that for every
form of level $1$ we have 
\[ \slope(f) < \frac{k}{4}.\]
For integral slopes (as we point out below, the slopes seem to be almost
always integral), this is in fact equivalent to our inequality 
\[ \slope(f) \leq \frac{k-1}{4}.\]
For $p=2$, Kevin Buzzard has formulated a conjectural description of
\emph{all} the slopes that implies, in particular, that the inequality
\[ \slope(f) \leq \frac{k-1}{3}\]
always holds for forms of level $1$, i.e., that there are no exceptional
forms for $p=2$. (Buzzard's conjectural description also implies that all
the $2$-adic slopes for level $1$ are integral.)

These hints suggest that something is going on. It seems to us that the
fact that the inequality is so often true demands some explanation. In
particular, one would like to know whether there is something special about
the cases where it fails.

\begin{question}
  Can one identify a specific property of the modular forms or the Galois
  representations corresponding to exceptional pairs $(p,k)$?
\end{question}

The list of exceptional cases is itself suggestive. Consider the case of
$p=59$. For $k=16$, one finds a form of slope $1$, which is therefore an
exception to the inequality. There is no clear explanation for this first
exceptional slope, but the fact that the other counterexamples occur at
weights $k=46,76,106,136,\dots$ suggests a systematic pattern. A possible
interpretation of this pattern via the theory of $\Theta$-cycles developed
in \cite{naomicong} was suggested by Kevin Buzzard and will be discussed in
the next section.

Finally, we note that \emph{almost all the slopes we obtained are
  integers}. Of all the observations we make, this is the one that is most
likely to be merely an effect of the fact that we work only with small
primes. The location of the exceptions, however, suggests that something
else may be going on. Specifically, non-integral slopes occur in our
computations only for $p=59$ and $p=79$, the same primes for which
exceptional slopes occur. Furthermore, the fractional slopes we observe are
connected to the exceptional slopes, in the following remarkable way:

\begin{enumerate}
\item A weight $k$ for which there exists exceptional form with slope equal
  to $2$ is preceded by a weight $k-2$ for which the slope sequence
  contains two slopes equal to $1/2$. For $p = 59$, this happens for the
  pairs of weights $(74,76)$ and $(104,106)$; for $p=79$, weights
  $(116,118)$, $(122,124)$.
\item A weight $k$ for which there exists exceptional form with slope equal
  to $3$ is preceded by a weight $k-2$ for which the slope sequence
  contains two slopes equal to $3/2$, and that weight is preceded by a
  weight $k-4$ whose slope sequence contains two slopes equal to $1/2$.
  This happens for $p=59$ and weights $(132,134,136)$ and $(162,164,166)$;
  for $p=79$, weights $(194,196,198)$ and $(200,202,204)$.
\item A weight $k$ for which there exists exceptional form with slope equal
  to $4$ is preceded by a weight $k-2$ for which the slope sequence
  contains two slopes equal to $5/2$, that weight is preceded by a weight
  $k-4$ whose slope sequence contains two slopes equal to $3/2$, and that
  weight is preceded by a weight $k-6$ whose slope sequence contains two
  slopes equal to $1/2$. This happens for $p=59$ and weights
  $(190,192,194,196)$ and $(220,222,224,226)$; for $p=79$ and weights
  $(272,274,276,278)$ and $(278,280,282,284)$.
\item And so on. An exceptional form of slope $n$ is and weight $k$ is
  accompanied by a ``train'' of pairs of forms of slope 
  \[ \frac{2n-3}{2},\: \frac{2n-5}{2},\: \frac{2n-7}{2},\: \dots,\: \frac12\]
  and weight
  \[ k-2,\: k-4,\: k-6,\: \dots,\: k-2(n-1)\]
  (one pair for each weight).  See the tables in section~\ref{tables} for
  the actual slope sequences.
\end{enumerate}

These patterns can overlap without interfering. For example, there are
exceptional forms of slope $6$ for $p=79$ and weights $438$ and $444$. Each
has its trail of weights for which fractional slopes exist. Because of the
exceptional form at weight $438$, there are forms of slope $1/2$ at weight
$428$ and at each subsequent weight, up to forms of slope $9/2$ at weight
$436$. Because of the exceptional form at weight $444$, there are forms of
slope $1/2$ at weight $434$ and at each subsequent weight, up to forms of
slope $9/2$ at weight $442$. Hence, for example, at weight $436$ we have
\emph{both} a pair of forms of slope $3/2$ and a pair of forms of slope
$9/2$.

This suggests a final question:

\begin{question}
  How often are the slopes integral? What is the connection between
  non-integral slopes and exceptional slopes?
\end{question}

\section{$\Theta$-cycles?}

In this section we look a little more closely at the exceptional cases that
occur for $p=59$ and $p=79$. The relevant slope sequences for $p=59$ and
for $p=79$ are listed in the tables in section~\ref{tables} The exceptional
and the fractional slopes are printed in boldface.

The list of weights at which exceptional slopes occur strongly suggest a
connection to the theory of $\Theta$-cycles, as in \cite{naomicong}. We
recall the basic ideas. As above, we restrict to the case of level $N=1$.
The $\Theta$ operator is the operator that acts on $q$-expansions as
$q\frac{d}{dq}$, so that
\[ \Theta\left( \sum a_nq^n\right) = \sum na_nq^n.\]
(If we think of modular forms as functions on the complex upper half-plane,
then $q=e^{2\pi i z}$ and $\Theta$ is just $\frac{d}{dz}$.) As is well
known (see \cite{story}), if $f$ is a modular form then $\Theta
f$ is \emph{not} a modular form, though it is ``almost'' modular in some
sense (the $p$-adic story is a little different: see \cite{ka350,1304,E2}).
On the other hand, $\Theta$ does define an operator on modular forms modulo
$p$, in which case it maps forms of weight $k$ to forms of weight $k+p+1$
(see \cite{naomicong,karesult}).

When one considers modular forms modulo $p$, it is possible for forms whose
weight differs by a multiple of $p-1$ to have identical $q$-expansions.
(Basically, this is because the $q$-expansion of the Eisenstein series
$E_{p-1}$ is congruent to $1$ modulo $p$.) Thus, if $f$ is a modular form,
we define its \emph{filtration} $w(f)$ to be the minimal weight for which
there is a modular form whose $q$-expansion is congruent modulo $p$ to the
$q$-expansion of $f$. The theory of $\Theta$-cycles describes what happens
to the filtration under the $\Theta$ operator.

The basic facts are the following. Let $f$ be a modular form of weight $k$
on $\Gamma_0(N)$. If $f$ is not ordinary (i.e., $U(f)\equiv 0 \pmod p$,
which is the case in all the examples we will consider), then
$\Theta^{p-1}f\equiv f\pmod p$, and hence $w(\Theta^{p-1}f)=w(f)$. This is
why one speaks of a ``$\Theta$-\emph{cycle}.''  The filtration normally
increases by $p+1$ each time we compute $\Theta$; specifically, if $w(f)$
is not divisible by $p$, then $w(\Theta f)=w(f)+p+1$. If, on the other
hand, $w(f)$ is divisible by $p$, the filtration goes down, so that
$w(\Theta f)=w(f)+p+1-n(p-1)$, where $n$ is some integer. What values $n$
can assume is completely determined in \cite{naomicong}. The case that is
relevant here is the one where $w(f)=k$, $4\leq k \leq p-1$, and $f$ is not
ordinary.  In this case, the cycle looks as follows. First,
\[ w(f)=k, w(\Theta f)=k+p+1, \dots,
  w(\Theta^{p-k}f)=p^2-(k-1)p.\]
Then, since this is divisible by $p$, the filtration falls, with $n=p+2-k$,
so 
\[ w(\Theta^{p-k+1}f)=p+3-k, \dots, w(\Theta^{p-2}f)=(k-2)p,\]
and finally $w(\Theta^{p-1}f)=k$, closing the cycle.

Notice that this whole theory refers \emph{only} to modular forms modulo
$p$. If $f$ is an eigenform, then all of the $\Theta^i f$ are eigenforms
modulo $p$; by the Deligne-Serre Lemma (see \cite[Lemma 6.11]{DS} or
\cite[Prop. 1.2.2]{glenn}), they lift to eigenforms in characteristic zero
(but recall that the $\Theta^i f$ are not themselves modular forms, so the
lifts will only be congruent to them). We know these eigenforms will have
positive slope (a property which is ``visible'' modulo $p$), but there
seems to no reason to predict anything further about their slope.

Consider now the exceptional slopes for $p=59$. For the case of
weight $16$, the occurrence of a large slope seems to be ``accidental,'' but
its ``propagation'' to higher weights seems to be linked a
$\Theta$-cycle. Let $f_0$ be the (unique) form of weight $16$ and level
$1$; its $59$-adic slope is $1$, because its $59$-th Fourier coefficient
is divisible (once) by $59$. For each $i$, let $f_i$ denote the minimal
weight modular form which is congruent modulo $p$ to $\Theta^i f$.
Then we have:

\begin{itemize}
\item $f_1$ is of weight $76$,
\item $f_2$ is of weight $136$,
\item $f_3$ is of weight $196$,
\item \dots
\item $f_{43}$ is of weight $2596=22\times59$,
\item $f_{44}$ is of weight $46$,
\item $f_{45}$ is of weight $106$,
\item $f_{46}$ is of weight $166$,
\item \dots
\item $f_{57}$ is of weight $826=14\times 59$,
\item finally, $f_{58}=f_0$.
\end{itemize}

Note that the smaller weights (up to $496$) on this list are precisely the
weights for which we found exceptional slopes.

Because the reductions modulo $p$ of these forms are in the image of
$\Theta$, all of the liftings of forms in the $\Theta$-cycle must have
positive slope, but the lifting theory suggests no reason to expect the
slope to increase along the cycle. Nor, in fact, is there any reason to
expect the slopes to be integers. What we observe, however, is the
following:

\begin{itemize}
\item The exceptional form of weight $76$ is congruent modulo $59$ (but not
  modulo $59^2$) to $\Theta f_0$ \emph{and has slope $2$}.
\item The exceptional form of weight $136$ is congruent modulo $59$ (but
  not modulo $59^2$) to $\Theta^2 f_0$ \emph{and has slope $3$}.
\item \dots
\item The exceptional form of weight $46$ is congruent modulo $59$ (but not
  modulo $59^2$) to $\Theta^{44} f_0$ \emph{and has slope $1$}.
\item The exceptional form of weight $106$ is congruent modulo $59$ (but
  not modulo $59^2$) to $\Theta^{45} f_0$ \emph{and has slope $2$}.
\item \dots
\end{itemize}

In other words, the forms we get are only connected via $\Theta$ modulo
$59$, but their slopes increase by $1$ as one might expect under $\Theta$,
except that at the point where the filtration falls back down to $46$ the
slope falls back down to $1$. If the pattern continues to hold for the full
cycle, every one of the resulting forms will be exceptional (i.e., their
supersingularities will be bigger than $1/60$), though their
supersingularities get closer and closer to $1/60$ as the weight gets
larger.

Also interesting is to ask whether one still gets exceptional slopes beyond
the $\Theta$-cycle; for example, is there an exceptional form for $p=59$ of
weight $886$? [Computation still in progress!]

Similarly, we can examine the exceptional slopes for $p=79$. Again, the
exceptional slopes seem to fit into a $\Theta$-cycle starting from a form
of weight $48$ (the weights in the cycle would be $48$, $128$, $208$,
\dots, $2528$, $34$, $114$, $194$, \dots, $3634$, $48$). Again, all we see
is the lower end of the cycle, and again, if the pattern persists
throughout the cycle every one of these forms would be exceptional.

If the existence of such ``$\Theta$-cycles'' of slopes reflects a general
phenomenon, then this would suggest that the set of possible slopes for a
given $p$ and varying $k$ has far more structure than predicted, say, in
\cite{conj}. This agrees with Kevin Buzzard's conjectures regarding the
case $p=2$.

\section{Tables of Exceptional and Fractional Slopes}\label{tables}

The tables that follow give the lower halves of the slope sequences for
$p=59$ and $p=79$ and all the weights for which there is either a
fractional or an exceptional slope. To make it easier for the reader, we
have printed the fractional and exceptional slopes in bold. The full
tables, containing all the slopes sequences for $p\leq 100$ and $k\leq
500$, can be found on the web at
\texttt{http://www.colby.edu/personal/fqgouvea/slopes/}. 

First the table for $p=59$:

\begin{longtable}{cl}
\toprule
Weight $k$ & Slope Sequence for $p=59$ (lower half) \\ 
\midrule
\endhead
\bottomrule
\endfoot
16         & (\textbf{1}) \\
46         & (0, 0, \textbf{1}) \\
74         & (\textbf{1/2}, \textbf{1/2}, 1, 1, 1)\\
76         & (0, 1, 1, 1, 1, \textbf{2}) \\
104        & (0, 0, \textbf{1/2}, \textbf{1/2}, 1, 1, 1, 1)\\
106        & (0, 0, 0, 0, 1, 1, 1, \textbf{2}) \\
132        & (\textbf{1/2}, \textbf{1/2}, 1, 1, 1, 2, 2, 2, 2, 2, 2)\\
134        & (0, 1, 1, 1, 1, \textbf{3/2}, \textbf{3/2}, 2, 2, 2)\\
136        & (0, 1, 1, 1, 1, 1, 2, 2, 2, 2, \textbf{3}) \\ 
162        & (0, 0, \textbf{1/2}, \textbf{1/2}, 1, 1, 1, 1, 2, 2, 2, 2, 2) \\
164        & (0, 0, 0, 0, 1, 1, 1, \textbf{3/2}, \textbf{3/2}, 2, 2, 2, 2) \\
166        & (0, 0, 0, 1, 1, 1, 1, 1, 1, 2, 2, 2, \textbf{3}) \\
190        & (\textbf{1/2}, \textbf{1/2}, 1, 1, 1, 2, 2, 2, 2, 2, 2, 3, 3, 3, 3) \\
192        & (0, 1, 1, 1, 1, \textbf{3/2}, \textbf{3/2}, 2, 2, 2, 3, 3, 3, 3, 3, 3) \\
194        & (0, 1, 1, 1, 1, 1, 2, 2, 2, 2, \textbf{5/2}, \textbf{5/2}, 3, 3, 3) \\
196        & (0, 1, 1, 1, 1, 1, 2, 2, 2, 2, 2, 3, 3, 3, 3, \textbf{4}) \\
220        & (0, 0, \textbf{1/2}, \textbf{1/2}, 1, 1, 1, 1, 2, 2, 2, 2, 2, 3, 3,\\
           & 3, 3, 3)\\
222        & (0, 0, 0, 0, 1, 1, 1, \textbf{3/2}, \textbf{3/2}, 2, 2, 2, 2, 3, 3, \\
           & 3, 3, 3)\\
224        & (0, 0, 0, 1, 1, 1, 1, 1, 1, 2, 2, 2, \textbf{5/2}, \textbf{5/2}, 3, \\
           & 3, 3, 3)\\
226        & (0, 0, 0, 0, 1, 1, 1, 1, 2, 2, 2, 2, 2, 2, 3, 3, 3, \textbf{4})\\
248        & (\textbf{1/2}, \textbf{1/2}, 1, 1, 1, 2, 2, 2, 2, 2, 2, 3, 3, 3, 3, \\
           & 4, 4, 4, 4, 4)\\
250        & (0, 1, 1, 1, 1, \textbf{3/2}, \textbf{3/2}, 2, 2, 2, 3, 3, 3, 3, 3, \\
           & 3, 4, 4, 4, 4)\\
252        & (0, 1, 1, 1, 1, 1, 2, 2, 2, 2, \textbf{5/2}, \textbf{5/2}, 3, 3, 3, \\
           & 4, 4, 4, 4, 4, 4)\\
254        & (0, 1, 1, 1, 1, 1, 2, 2, 2, 2, 2, 3, 3, 3, 3, \textbf{7/2}, \textbf{7/2}, \\
           & 4, 4, 4)\\
256        & (0, 0, 1, 1, 1, 1, 2, 2, 2, 2, 2, 3, 3, 3, 3, 3, 4,\\
           & 4, 4, 4, \textbf{5})\\
278        & (0, 0, \textbf{1/2}, \textbf{1/2}, 1, 1, 1, 1, 2, 2, 2, 2, 2, 3, 3, \\
           & 3, 3, 3, 4, 4, 4, 4)\\
280        & (0, 0, 0, 0, 1, 1, 1, \textbf{3/2}, \textbf{3/2}, 2, 2, 2, 2, 3, 3, \\
           & 3, 3, 3, 4, 4, 4, 4, 4)\\
282        & (0, 0, 0, 1, 1, 1, 1, 1, 1, 2, 2, 2, \textbf{5/2}, \textbf{5/2}, 3, \\
           & 3, 3, 3, 4, 4, 4, 4, 4)\\
284        & (0, 0, 0, 0, 1, 1, 1, 1, 2, 2, 2, 2, 2, 2, 3, 3, 3, \textbf{7/2}, \textbf{7/2}, \\
           & 4, 4, 4, 4)\\
286        & (0, 0, 0, 0, 1, 1, 1, 1, 1, 2, 2, 2, 2, 3, 3, 3, 3, 3, 3, 4, 4, 4, \textbf{5})\\
306        & (\textbf{1/2}, \textbf{1/2}, 1, 1, 1, 2, 2, 2, 2, 2, 2, 3, 3, 3, 3, \\
           & 4, 4, 4, 4, 4, 5, 5, 5, 5, 5)\\
308        & (0, 1, 1, 1, 1, \textbf{3/2}, \textbf{3/2}, 2, 2, 2, 3, 3, 3, 3, 3, \\
           & 3, 4, 4, 4, 4, 5, 5, 5, 5, 5)\\
310        & (0, 1, 1, 1, 1, 1, 2, 2, 2, 2, \textbf{5/2}, \textbf{5/2}, 3, 3, 3, \\
           & 4, 4, 4, 4, 4, 4, 5, 5, 5, 5)\\
312        & (0, 1, 1, 1, 1, 1, 2, 2, 2, 2, 2, 3, 3, 3, 3, \textbf{7/2}, \textbf{7/2}, \\
           & 4, 4, 4, 5, 5, 5, 5, 5, 5)\\
314        & (0, 0, 1, 1, 1, 1, 2, 2, 2, 2, 2, 3, 3, 3, 3, 3, 4, 4, 4, 4, \textbf{9/2}, \\
           & \textbf{9/2}, 5, 5, 5)\\
316        & (0, 1, 1, 1, 1, 1, 1, 2, 2, 2, 2, 3, 3, 3, 3, 3, 4, 4, 4, 4, 4, 5, 5, \\
           & 5, 5, \textbf{6})\\
336        & (0, 0, \textbf{1/2}, \textbf{1/2}, 1, 1, 1, 1, 2, 2, 2, 2, 2, 3, 3, 3, \\
           & 3, 3, 4, 4, 4, 4, 5, 5, 5, 5, 5, 5)\\
338        & (0, 0, 0, 0, 1, 1, 1, \textbf{3/2}, \textbf{3/2}, 2, 2, 2, 2, 3, 3, 3, \\
           & 3, 3, 4, 4, 4, 4, 4, 5, 5, 5, 5)\\
340        & (0, 0, 0, 1, 1, 1, 1, 1, 1, 2, 2, 2, \textbf{5/2}, \textbf{5/2}, 3, 3, \\
           & 3, 3, 4, 4, 4, 4, 4, 5, 5, 5, 5, 5)\\
342        & (0, 0, 0, 0, 1, 1, 1, 1, 2, 2, 2, 2, 2, 2, 3, 3, 3, \textbf{7/2}, \textbf{7/2}, \\
           & 4, 4, 4, 4, 5, 5, 5, 5, 5)\\
344        & (0, 0, 0, 0, 1, 1, 1, 1, 1, 2, 2, 2, 2, 3, 3, 3, 3, 3, 3, 4, 4, 4, \textbf{9/2}, \\
           & \textbf{9/2}, 5, 5, 5, 5)\\
346        & (0, 0, 0, 0, 1, 1, 1, 1, 1, 2, 2, 2, 2, 2, 3, 3, 3, 3, 4, 4, 4, 4, 4, \\
           & 4, 5, 5, 5, \textbf{6})\\
364        & (\textbf{1/2}, \textbf{1/2}, 1, 1, 1, 2, 2, 2, 2, 2, 2, 3, 3, 3, 3, 4, \\
           & 4, 4, 4, 4, 5, 5, 5, 5, 5, 6, 6, 6, 6, 6)\\
366        & (0, 1, 1, 1, 1, \textbf{3/2}, \textbf{3/2}, 2, 2, 2, 3, 3, 3, 3, 3, 3, \\
           & 4, 4, 4, 4, 5, 5, 5, 5, 5, 6, 6, 6, 6, 6)\\
368        & (0, 1, 1, 1, 1, 1, 2, 2, 2, 2, \textbf{5/2}, \textbf{5/2}, 3, 3, 3, 4, \\
           & 4, 4, 4, 4, 4, 5, 5, 5, 5, 6, 6, 6, 6, 6)\\
370        & (0, 1, 1, 1, 1, 1, 2, 2, 2, 2, 2, 3, 3, 3, 3, \textbf{7/2}, \textbf{7/2}, \\
           & 4, 4, 4, 5, 5, 5, 5, 5, 5, 6, 6, 6, 6)\\
372        & (0, 0, 1, 1, 1, 1, 2, 2, 2, 2, 2, 3, 3, 3, 3, 3, 4, 4, 4, 4, \textbf{9/2}, \\
           & \textbf{9/2}, 5, 5, 5, 6, 6, 6, 6, 6, 6)\\
374        & (0, 1, 1, 1, 1, 1, 1, 2, 2, 2, 2, 3, 3, 3, 3, 3, 4, 4, 4, 4, 4, 5, 5, \\
           & 5, 5, \textbf{11/2}, \textbf{11/2}, 6, 6, 6)\\
376        & (0, 0, 1, 1, 1, 1, 2, 2, 2, 2, 2, 2, 3, 3, 3, 3, 4, 4, 4, 4, 4, 5, 5, \\
           & 5, 5, 5, 6, 6, 6, 6, \textbf{7})\\
394        & (0, 0, \textbf{1/2}, \textbf{1/2}, 1, 1, 1, 1, 2, 2, 2, 2, 2, 3, 3, 3, \\
           & 3, 3, 4, 4, 4, 4, 5, 5, 5, 5, 5, 5, 6, 6, 6, 6)\\
396        & (0, 0, 0, 0, 1, 1, 1, \textbf{3/2}, \textbf{3/2}, 2, 2, 2, 2, 3, 3, 3, \\
           & 3, 3, 4, 4, 4, 4, 4, 5, 5, 5, 5, 6, 6, 6, 6, 6, 6)\\
398        & (0, 0, 0, 1, 1, 1, 1, 1, 1, 2, 2, 2, \textbf{5/2}, \textbf{5/2}, 3, 3, \\
           & 3, 3, 4, 4, 4, 4, 4, 5, 5, 5, 5, 5, 6, 6, 6, 6)\\
400        & (0, 0, 0, 0, 1, 1, 1, 1, 2, 2, 2, 2, 2, 2, 3, 3, 3, \textbf{7/2}, \textbf{7/2}, \\
           & 4, 4, 4, 4, 5, 5, 5, 5, 5, 6, 6, 6, 6, 6)\\
402        & (0, 0, 0, 0, 1, 1, 1, 1, 1, 2, 2, 2, 2, 3, 3, 3, 3, 3, 3, 4, 4, 4, \textbf{9/2}, \\
           & \textbf{9/2}, 5, 5, 5, 5, 6, 6, 6, 6, 6)\\
404        & (0, 0, 0, 0, 1, 1, 1, 1, 1, 2, 2, 2, 2, 2, 3, 3, 3, 3, 4, 4, 4, 4, 4, \\
           & 4, 5, 5, 5, \textbf{11/2}, \textbf{11/2}, 6, 6, 6, 6)\\
406        & (0, 0, 0, 0, 1, 1, 1, 1, 1, 2, 2, 2, 2, 2, 3, 3, 3, 3, 3, 4, 4, 4, 4, \\
           & 5, 5, 5, 5, 5, 5, 6, 6, 6, \textbf{7})\\
422        & (\textbf{1/2}, \textbf{1/2}, 1, 1, 1, 2, 2, 2, 2, 2, 2, 3, 3, 3, 3, 4, \\
           & 4, 4, 4, 4, 5, 5, 5, 5, 5, 6, 6, 6, 6, 6, 7, 7, 7, 7)\\
424        & (0, 1, 1, 1, 1, \textbf{3/2}, \textbf{3/2}, 2, 2, 2, 3, 3, 3, 3, 3, 3, \\
           & 4, 4, 4, 4, 5, 5, 5, 5, 5, 6, 6, 6, 6, 6, 7, 7, 7, 7, 7)\\
426        & (0, 1, 1, 1, 1, 1, 2, 2, 2, 2, \textbf{5/2}, \textbf{5/2}, 3, 3, 3, 4, \\
           & 4, 4, 4, 4, 4, 5, 5, 5, 5, 6, 6, 6, 6, 6, 7, 7, 7, 7, 7)\\
428        & (0, 1, 1, 1, 1, 1, 2, 2, 2, 2, 2, 3, 3, 3, 3, \textbf{7/2}, \textbf{7/2}, \\
           & 4, 4, 4, 5, 5, 5, 5, 5, 5, 6, 6, 6, 6, 7, 7, 7, 7, 7)\\
430        & (0, 0, 1, 1, 1, 1, 2, 2, 2, 2, 2, 3, 3, 3, 3, 3, 4, 4, 4, 4, \textbf{9/2}, \\
           & \textbf{9/2}, 5, 5, 5, 6, 6, 6, 6, 6, 6, 7, 7, 7, 7)\\
432        & (0, 1, 1, 1, 1, 1, 1, 2, 2, 2, 2, 3, 3, 3, 3, 3, 4, 4, 4, 4, 4, 5, 5, \\
           & 5, 5, \textbf{11/2}, \textbf{11/2}, 6, 6, 6, 7, 7, 7, 7, 7, 7)\\
434        & (0, 0, 1, 1, 1, 1, 2, 2, 2, 2, 2, 2, 3, 3, 3, 3, 4, 4, 4, 4, 4, 5, 5, \\
           & 5, 5, 5, 6, 6, 6, 6, \textbf{13/2}, \textbf{13/2}, 7, 7, 7)\\
436        & (0, 0, 1, 1, 1, 1, 1, 2, 2, 2, 2, 3, 3, 3, 3, 3, 3, 4, 4, 4, 4, 5, 5, \\
           & 5, 5, 5, 6, 6, 6, 6, 6, 7, 7, 7, 7, \textbf{8})\\
452       & (0, 0, \textbf{1/2}, \textbf{1/2}, 1, 1, 1, 1, 2, 2, 2, 2, 2, 3, 3, 3, \\
           & 3, 3, 4, 4, 4, 4, 5, 5, 5, 5, 5, 5, 6, 6, 6, 6, 7, 7, 7, 7, 7)\\
454       & (0, 0, 0, 0, 1, 1, 1, \textbf{3/2}, \textbf{3/2}, 2, 2, 2, 2, 3, 3, 3, \\
           & 3, 3, 4, 4, 4, 4, 4, 5, 5, 5, 5, 6, 6, 6, 6, 6, 6, 7, 7, 7, 7)\\
456       & (0, 0, 0, 1, 1, 1, 1, 1, 1, 2, 2, 2, \textbf{5/2}, \textbf{5/2}, 3, 3, \\
           & 3, 3, 4, 4, 4, 4, 4, 5, 5, 5, 5, 5, 6, 6, 6, 6, 7, 7, 7, 7, 7, 7)\\
458       & (0, 0, 0, 0, 1, 1, 1, 1, 2, 2, 2, 2, 2, 2, 3, 3, 3, \textbf{7/2}, \textbf{7/2}, \\
           & 4, 4, 4, 4, 5, 5, 5, 5, 5, 6, 6, 6, 6, 6, 7, 7, 7, 7)\\
460       & (0, 0, 0, 0, 1, 1, 1, 1, 1, 2, 2, 2, 2, 3, 3, 3, 3, 3, 3, 4, 4, 4, \textbf{9/2}, \\
           & \textbf{9/2}, 5, 5, 5, 5, 6, 6, 6, 6, 6, 7, 7, 7, 7, 7)\\
462       & (0, 0, 0, 0, 1, 1, 1, 1, 1, 2, 2, 2, 2, 2, 3, 3, 3, 3, 4, 4, 4, 4, 4, \\
           & 4, 5, 5, 5, \textbf{11/2}, \textbf{11/2}, 6, 6, 6, 6, 7, 7, 7, 7, 7)\\
464       & (0, 0, 0, 0, 1, 1, 1, 1, 1, 2, 2, 2, 2, 2, 3, 3, 3, 3, 3, 4, 4, 4, 4, \\
           & 5, 5, 5, 5, 5, 5, 6, 6, 6, \textbf{13/2}, \textbf{13/2}, 7, 7, 7, 7)\\
466       & (0, 0, 0, 0, 0, 1, 1, 1, 1, 2, 2, 2, 2, 2, 3, 3, 3, 3, 3, 4, 4, 4, 4, \\
           & 4, 5, 5, 5, 5, 6, 6, 6, 6, 6, 6, 7, 7, 7, \textbf{8})\\
480       & (\textbf{1/2}, \textbf{1/2}, 1, 1, 1, 2, 2, 2, 2, 2, 2, 3, 3, 3, 3, 4, \\
           & 4, 4, 4, 4, 5, 5, 5, 5, 5, 6, 6, 6, 6, 6, 7, 7, 7, 7, 7, 7, 7, 7, 7, 7)\\
482       & (0, 1, 1, 1, 1, \textbf{3/2}, \textbf{3/2}, 2, 2, 2, 3, 3, 3, 3, 3, 3, \\
           & 4, 4, 4, 4, 5, 5, 5, 5, 5, 6, 6, 6, 6, 6, 7, 7, 7, 7, 7, 8, 8, 8, 8)\\
484       & (0, 1, 1, 1, 1, 1, 2, 2, 2, 2, \textbf{5/2}, \textbf{5/2}, 3, 3, 3, 4, \\
           & 4, 4, 4, 4, 4, 5, 5, 5, 5, 6, 6, 6, 6, 6, 7, 7, 7, 7, 7, 8, 8, 8, 8, 8)\\
486       & (0, 1, 1, 1, 1, 1, 2, 2, 2, 2, 2, 3, 3, 3, 3, \textbf{7/2}, \textbf{7/2}, \\
           & 4, 4, 4, 5, 5, 5, 5, 5, 5, 6, 6, 6, 6, 7, 7, 7, 7, 7, 8, 8, 8, 8, 8)\\
488       & (0, 0, 1, 1, 1, 1, 2, 2, 2, 2, 2, 3, 3, 3, 3, 3, 4, 4, 4, 4, \textbf{9/2}, \\
           & \textbf{9/2}, 5, 5, 5, 6, 6, 6, 6, 6, 6, 7, 7, 7, 7, 8, 8, 8, 8, 8)\\
490       & (0, 1, 1, 1, 1, 1, 1, 2, 2, 2, 2, 3, 3, 3, 3, 3, 4, 4, 4, 4, 4, 5, 5, \\
           & 5, 5, \textbf{11/2}, \textbf{11/2}, 6, 6, 6, 7, 7, 7, 7, 7, 7, 8, 8, 8, 8)\\
492       & (0, 0, 1, 1, 1, 1, 2, 2, 2, 2, 2, 2, 3, 3, 3, 3, 4, 4, 4, 4, 4, 5, 5, \\
           & 5, 5, 5, 6, 6, 6, 6, \textbf{13/2}, \textbf{13/2}, 7, 7, 7, 8, 8, 8, 8, 8, 8)\\
494       & (0, 0, 1, 1, 1, 1, 1, 2, 2, 2, 2, 3, 3, 3, 3, 3, 3, 4, 4, 4, 4, 5, 5, \\
           & 5, 5, 5, 6, 6, 6, 6, 6, 7, 7, 7, 7, \textbf{15/2}, \textbf{15/2}, 8, 8, 8)\\
496       & (0, 0, 1, 1, 1, 1, 1, 2, 2, 2, 2, 2, 3, 3, 3, 3, 4, 4, 4, 4, 4, 4, 5, \\
           & 5, 5, 5, 6, 6, 6, 6, 6, 7, 7, 7, 7, 7, 8, 8, 8, 8, \textbf{9})
\end{longtable}

Now the table for $p=79$:

\begin{longtable}{cll}
\toprule
Weight $k$ & Slope sequence for $p=79$(lower half)\\ 
\midrule
\endhead
\bottomrule
\endfoot
38   &  (0, \textbf{1})\\
44   &  (0, 0, \textbf{1})\\
116  &  (0, \textbf{1/2}, \textbf{1/2}, 1, 1, 1, 1, 1, 1)\\
118  &  (0, 0, 0, 1, 1, 1, 1, 1, \textbf{2})\\
122  &  (0, 0, \textbf{1/2}, \textbf{1/2}, 1, 1, 1, 1, 1)\\
124  &  (0, 0, 0, 1, 1, 1, 1, 1, 1, \textbf{2})\\
194  &  (0, \textbf{1/2}, \textbf{1/2}, 1, 1, 1, 1, 1, 1, 2, 2, 2, 2, 2, 2)\\
196  &  (0, 0, 0, 1, 1, 1, 1, 1, \textbf{3/2}, \textbf{3/2}, 2, 2, 2, 2, 2, 2)\\
198  &  (0, 0, 0, 1, 1, 1, 1, 1, 1, 1, 2, 2, 2, 2, 2, \textbf{3})\\
200  &  (0, 0, \textbf{1/2}, \textbf{1/2}, 1, 1, 1, 1, 1, 2, 2, 2, 2, 2, 2, 2)\\
202  &  (0, 0, 0, 1, 1, 1, 1, 1, 1, \textbf{3/2}, \textbf{3/2}, 2, 2, 2, 2, 2)\\
204  &  (0, 0, 0, 0, 1, 1, 1, 1, 1, 1, 2, 2, 2, 2, 2, 2, \textbf{3})\\
272  &  (0, \textbf{1/2}, \textbf{1/2}, 1, 1, 1, 1, 1, 1, 2, 2, 2, 2, 2, 2, 3, 3, 3,\\
     &    3, 3, 3, 3)\\
274  &  (0, 0, 0, 1, 1, 1, 1, 1, \textbf{3/2}, \textbf{3/2}, 2, 2, 2, 2, 2, 2, 3, 3,\\
     &    3, 3, 3, 3)\\
276  &  (0, 0, 0, 1, 1, 1, 1, 1, 1, 1, 2, 2, 2, 2, 2, \textbf{5/2}, \textbf{5/2}, 3,\\
     &    3, 3, 3, 3, 3)\\
278  &  (0, 0, \textbf{1/2}, \textbf{1/2}, 1, 1, 1, 1, 1, 2, 2, 2, 2, 2, 2, 2, 3, 3,\\
     &    3, 3, 3, \textbf{4})\\
280  &  (0, 0, 0, 1, 1, 1, 1, 1, 1, \textbf{3/2}, \textbf{3/2}, 2, 2, 2, 2, 2, 3, 3,\\
     &    3, 3, 3, 3, 3)\\
282  &  (0, 0, 0, 0, 1, 1, 1, 1, 1, 1, 2, 2, 2, 2, 2, 2, \textbf{5/2}, \textbf{5/2},\\
     &    3, 3, 3, 3, 3)\\
284  &  (0, 0, 0, 1, 1, 1, 1, 1, 1, 1, 2, 2, 2, 2, 2, 2, 3, 3, 3, 3,\\
     &    3, 3, \textbf{4})\\
350  &  (0, \textbf{1/2}, \textbf{1/2}, 1, 1, 1, 1, 1, 1, 2, 2, 2, 2, 2, 2, 3, 3, 3,\\
     &    3, 3, 3, 3, 4, 4, 4, 4, 4, 4)\\
352  &  (0, 0, 0, 1, 1, 1, 1, 1, \textbf{3/2}, \textbf{3/2}, 2, 2, 2, 2, 2, 2, 3, 3, 3, 3,\\
     &   3, 3, 4, 4, 4, 4, 4, 4, 4)\\
354  &  (0, 0, 0, 1, 1, 1, 1, 1, 1, 1, 2, 2, 2, 2, 2, \textbf{5/2}, \textbf{5/2}, 3, 3, 3,\\
     &   3, 3, 3, 4, 4, 4, 4, 4, 4)\\
356  &  (0, 0, \textbf{1/2}, \textbf{1/2}, 1, 1, 1, 1, 1, 2, 2, 2, 2, 2, 2, 2, 3, 3, 3, 3,\\
     &   3, \textbf{7/2}, \textbf{7/2}, 4, 4, 4, 4, 4, 4)\\
358  &  (0, 0, 0, 1, 1, 1, 1, 1, 1, \textbf{3/2}, \textbf{3/2}, 2, 2, 2, 2, 2, 3, 3, 3, 3,\\
     &   3, 3, 3, 4, 4, 4, 4, 4, \textbf{5})\\
360  &  (0, 0, 0, 0, 1, 1, 1, 1, 1, 1, 2, 2, 2, 2, 2, 2, \textbf{5/2}, \textbf{5/2}, 3, 3,\\
     &   3, 3, 3, 4, 4, 4, 4, 4, 4, 4)\\
362  &  (0, 0, 0, 1, 1, 1, 1, 1, 1, 1, 2, 2, 2, 2, 2, 2, 3, 3, 3, 3, 3, 3,\\
     &   \textbf{7/2}, \textbf{7/2}, 4, 4, 4, 4, 4)\\
364  &  (0, 0, 0, 0, 1, 1, 1, 1, 1, 1, 2, 2, 2, 2, 2, 2, 2, 3, 3, 3, 3, 3,\\
     &   3, 4, 4, 4, 4, 4, 4, \textbf{5})\\
428  &  (0, \textbf{1/2}, \textbf{1/2}, 1, 1, 1, 1, 1, 1, 2, 2, 2, 2, 2, 2, 3, 3, 3, 3, 3,\\
     &   3, 3, 4, 4, 4, 4, 4, 4, 5, 5, 5, 5, 5, 5, 5)\\
430  &  (0, 0, 0, 1, 1, 1, 1, 1, \textbf{3/2}, \textbf{3/2}, 2, 2, 2, 2, 2, 2, 3, 3, 3, 3,\\
     &   3, 3, 4, 4, 4, 4, 4, 4, 4, 5, 5, 5, 5, 5, 5)\\
432  &  (0, 0, 0, 1, 1, 1, 1, 1, 1, 1, 2, 2, 2, 2, 2, \textbf{5/2}, \textbf{5/2}, 3, 3, 3,\\
     &   3, 3, 3, 4, 4, 4, 4, 4, 4, 5, 5, 5, 5, 5, 5, 5)\\
434  &  (0, 0, \textbf{1/2}, \textbf{1/2}, 1, 1, 1, 1, 1, 2, 2, 2, 2, 2, 2, 2, 3, 3, 3, 3,\\
     &   3, \textbf{7/2}, \textbf{7/2}, 4, 4, 4, 4, 4, 4, 5, 5, 5, 5, 5, 5)\\
436  &  (0, 0, 0, 1, 1, 1, 1, 1, 1, \textbf{3/2}, \textbf{3/2}, 2, 2, 2, 2, 2, 3, 3, 3, 3,\\
     &   3, 3, 3, 4, 4, 4, 4, 4, \textbf{9/2}, \textbf{9/2}, 5, 5, 5, 5, 5, 5)\\
438  &  (0, 0, 0, 0, 1, 1, 1, 1, 1, 1, 2, 2, 2, 2, 2, 2, \textbf{5/2}, \textbf{5/2}, 3, 3,\\
     &   3, 3, 3, 4, 4, 4, 4, 4, 4, 4, 5, 5, 5, 5, 5, \textbf{6})\\
440  &  (0, 0, 0, 1, 1, 1, 1, 1, 1, 1, 2, 2, 2, 2, 2, 2, 3, 3, 3, 3, 3, 3,\\
     &   \textbf{7/2}, \textbf{7/2}, 4, 4, 4, 4, 4, 5, 5, 5, 5, 5, 5, 5)\\
442  &  (0, 0, 0, 0, 1, 1, 1, 1, 1, 1, 2, 2, 2, 2, 2, 2, 2, 3, 3, 3, 3, 3,\\
     &   3, 4, 4, 4, 4, 4, 4, \textbf{9/2}, \textbf{9/2}, 5, 5, 5, 5, 5)\\
444  &  (0, 0, 0, 0, 1, 1, 1, 1, 1, 1, 1, 2, 2, 2, 2, 2, 2, 3, 3, 3, 3, 3,\\
     &   3, 3, 4, 4, 4, 4, 4, 4, 5, 5, 5, 5, 5, 5, \textbf{6})\\
\end{longtable}

\providecommand{\bysame}{\leavevmode\hbox to3em{\hrulefill}\thinspace}

\end{document}